\documentclass[11pt]{amsart}
\usepackage{graphicx}

\newtheorem{theorem}{Theorem}[section]

\newtheorem{proposition}{Proposition}[section]

\theoremstyle{definition}

\newcommand{\bbar}{\left[ \begin{array}}
\newcommand{\ebar}{\end{array} \right] }
\newcommand{\bdm}{\begin{displaymath}}
\newcommand{\edm}{\end{displaymath}}
\newcommand{\beq}{\begin{equation}}
\newcommand{\beqa}{\begin{eqnarray}}
\newcommand{\beqas}{\begin{eqnarray*}}
\newcommand{\eeq}{\end{equation}}
\newcommand{\eeqa}{\end{eqnarray}}
\newcommand{\eeqas}{\end{eqnarray*}}
\newcommand{\dd}{\textup{d}}

\newcommand{\C}{\mathbb C}

\newcommand{\NN}{\bar{N}}

\newcommand{\Khat}{\widehat{K}}

 \newcommand{\real}{{\bf R}}

\newcommand{\Ubar}{\bar{U}}

\newcommand{\Kbar}{\bar{K}}

\newcommand{\shat}{\hat{\sigma}}

\newcommand{\kbar}{\bar{\mathfrak{k}}}
\newcommand{\pbar}{\bar{\mathfrak{p}}}
\newcommand{\ubar}{\bar{{\mathfrak{u}}}}

\newcommand{\pp}{{\mathfrak{p}}}

\newcommand{\tbar}{{\bar{\tau}}}

\newcommand{\VV}{\mathcal{V}}

\newcommand{\hh}{\mathcal{H}}

\newcommand{\cc}{{\bf C}}

\theoremstyle{remark}

\numberwithin{equation}{section}

\begin{document}

\title[Grassmann Geometries and Integrable Systems]{Grassmann Geometries and Integrable Systems}

% Remove or comment out any unused author tags.
% author one information
\author{David Brander}
\address{Department of Mathematics\\ Faculty of Science\\Kobe University\\1-1, Rokkodai, Nada-ku, Kobe 657-8501\\ Japan}
\email{brander@math.kobe-u.ac.jp}

%\thanks{}

% Use this \subjclass if you are using amsart version 2.0 (December 1999).
\subjclass[2000]{Primary 53C42, 53B25; Secondary 37J35, 53C35}
% Use this one if you are using an older version of amsart.
%\subjclass{}
\date{}

% at present the "communicated by" line appears only in ERA and PROC
\commby{}

\dedicatory{}

\begin{abstract}
We describe how the loop group maps corresponding to special submanifolds associated
to integrable systems may be thought of as certain Grassmann submanifolds of infinite
dimensional homogeneous spaces.  In general, the associated families of special submanifolds
are certain Grassmann submanifolds.  An example is given from 
the recent article \cite{reflective}.
\end{abstract}

\maketitle
%***************************************

%***************************************************

\section{Introduction}
This article discusses some of the ideas in the article \cite{reflective},
where solutions to a certain loop group problem were studied. The emphasis
here is on the geometric interpretation of the solutions, rather than the
techniques for producing solutions.
 
 In 1996, Ferus and Pedit \cite{feruspedit1996}
  defined an integrable system involving a 3-involution
loop group, solutions of which are isometric immersions between space forms of 
different non-zero sectional curvature. They modified the Adler-Kostant-Symes (AKS)
theory (described in \cite{burstallpedit}) to show how to produce many solutions
by solving commuting ODEs on a finite dimensional vector space.

The present author later studied this system in 
\cite{brander2} and \cite{branderrossman}:
 it had several interesting
properties, including a relationship with
pluriharmonic maps.  \\

\textbf{Goal here:}  generalize the system
to arbitrary commuting involutions of any Lie group and 
identify the associated special submanifolds.\\

\textbf{Results:} briefly, we obtained:
\begin{itemize}
\item
Generalizations, to all reflective submanifolds, of results concerning 
isometric immersions of space forms; 
\item
In case of previous results, new proofs;
\item
And other new special submanifolds as integrable systems. 
\end{itemize}

\subsection{Motivation}
Other special submanifolds that have been studied with loop groups,
(e.g. harmonic maps into symmetric spaces,
 CMC surfaces, special Lagrangian surfaces etc),
are associated to loop groups with only two involutions.  
Therefore, it seemed that a system in a loop group with three
involutions might have some interesting  properties peculiar to this situation.

One such property, studied in \cite{brander2}, is as follows: solutions to three 
\emph{distinct} problems are obtained from the \emph{same} loop group map,
by evaluating the map within different ranges of the loop parameter $\lambda$.
This amounts to a kind of Lawson correspondence between solutions of these
problems, and shows that the problems of obtaining complete immersions
are equivalent for the three cases. \\

\begin{table}[here] 
  \begin{tabular}{|c|c|c|}  \hline
Parameter range & Induced sectional curvature  & Target space \\ \hline
$\lambda \in i \real^*$  & $c_\lambda \in (-\infty,0)$ & $S^{m+k}$  \\
$\lambda \in  \real^*$  & $c_\lambda \in [-1,0)$ & $H^{m+k}_{k}$  \\
$\lambda \in S^1$  & $c_\lambda \in (-\infty,-1]$ & $H^{m+k}$ \\
\hline
  \end{tabular}
\label{table1}
\end{table}

The table shows three different constant curvature 
Riemannian submanifolds of three different
space forms obtained by evaluating the \emph{same} loop group
map for values of the spectral parameter in $\real$, $i \real$ and $S^1$ \cite{brander2}.\\

%*********************************************
\section{Special Submanifolds and Loop Groups}
We first present an outline of how certain special submanifolds are associated to maps
into loop groups.

\subsection{Moving Frame Method}  
The basic concept of the moving frame method is encapsulated as follows:
 \begin{itemize}
  \item
  Given  $f: M \to G/H$, an immersed submanifold of a homogeneous space.
  \item
    Lift, $F: M \to G$, a \textbf{frame} for $f$.
      \item \textbf{Idea:} Choose $F$ which is adapted in some
  way to the geometry of $f$. 
\end{itemize} 

%!!!!!!!!!!!!!!!!!!!!!!!!!!!!!!!!!!!!!!!!!!!!!!!!!!!!!!!!!!
%!!!!!!!!!!!!!!!!!!!!!!!!!!!!!!!!!!!!!!!!!!!!!!!!!!!!!!!!!!
\begin{figure}[here]
\includegraphics[width=4cm]{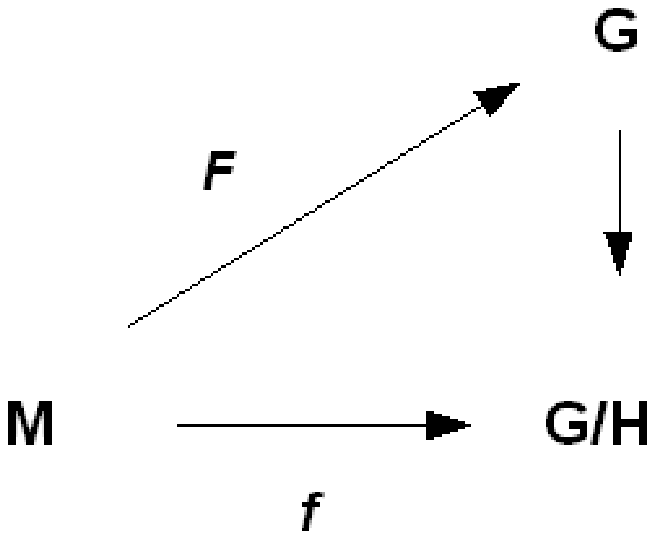}
\end{figure}

%-------------
\textbf{Example:} We illustrate this with a simple example.\\
{Special submanifold:} a flat immersion,
\bdm 
f: M = \real^2 \to S^3,
\edm
{Adapted frame:} $F: \real^2 \to SO(4)$,
\bdm
F := \bbar {cccc} e_1 & e_2 &  n & f \ebar, 
\edm
where  $e_i$ are an orthonormal basis for the tangent space to the immersion.

%--------
\subsection{The Maurer-Cartan Form}  
Given a frame $F: M \to G$,  for  $f : M \to G/H$,
the  \emph{Maurer-Cartan form}, $\alpha = F^{-1} \dd F \, \in \, \mathfrak{g} \otimes \Omega (M)$,
is the pull-back to $M$ of the Maurer-Cartan form of $G$.  
It is necessary that $\alpha$ satisfies the \emph{Maurer-Cartan equation}
     \beq \label{mceq}
    \dd \alpha + \alpha \wedge \alpha =0.
    \eeq
Conversely, if any $\alpha \,  \in \mathfrak{g} \otimes \Omega (M)$, satisfies
    (\ref{mceq}) then it is a basic fact from the theory of Lie groups that we
    can integrate $\alpha$ to obtain a map $F: M \to G$, whose Maurer-Cartan form
    is $\alpha$.  The map $F$ is determined up to an initial condition $F_0 \in G$. Changing
    this initial condition amounts to left multiplication by an element of $G$, which
    is to say an isometry of the homogeneous space $G/H$, and consequently we have the\\
    \textbf{Fundamental point:} $\alpha$ contains all geometric information about $f$.\\
  
\textbf{Example:} Returning to our previous example of flat surfaces in $S^3$, we
compute the Maurer-Cartan form of 
$F := \bbar {cccc} e_1 & e_2 &  n & f \ebar,$

\beqas
\alpha = F^{-1} \dd F &=& \bbar {c} e_1^T\\ e_2^T\\ n^T\\ f^T \ebar
 \cdot \bbar {cccc} \dd e_1 & \dd e_2 & \dd n & \dd f \ebar \\
&=& \bbar {ccc} \omega & \beta & \theta \\
          -\beta^t & 0 & 0 \\
          -\theta^t & 0 & 0 \ebar,
\eeqas
where the $2 \times 2$ matrix $\omega$ is the connection on the tangent bundle for $f$,
the $2 \times 1$ vector $\beta$ is the second fundamental form, and the $2 \times 1$ vector
$\theta$ is the coframe.

Computing the Maurer-Cartan equation
 $\dd \alpha + \alpha \wedge \alpha = 0$, the three components above give the following
 three equations:
\beqa
\dd \omega + \omega \wedge \omega -\beta \wedge \beta^t -\theta \wedge \theta ^t = 0, \label{eq1}\\
\dd \beta + \omega \wedge \beta = 0, \label{eq2}\\
\dd \theta + \omega \wedge \theta = 0. \label{eq3}
\eeqa
The assumption that the induced metric is flat is given by a further equation,
 \textbf{Flatness:}
\bdm
 \dd \omega + \omega \wedge \omega = 0.
\edm

\subsection{Parameterised Families of Frames}
Now suppose we introduce a complex parameter $\lambda$ in our example by setting:
\bdm 
\alpha_{\lambda} = \bbar  {ccc} \omega & {\lambda }\beta & {\lambda } \theta \\
          - {\lambda }\beta^t & 0 & 0 \\
          -{\lambda }\theta^t & 0 & 0 \ebar
    = a_0 + a_1 \lambda .
\edm
 Then 
$\dd \alpha_\lambda + \alpha_\lambda \wedge \alpha_\lambda = 0$   $\Leftrightarrow $ 
$\dd \omega + \omega \wedge \omega -{\lambda }^2(\beta \wedge \beta^t + \theta \wedge \theta ^t) = 0$, plus (\ref{eq2}) and (\ref{eq3}). It follows that we have the following 
equivalence:
\bdm
\dd \alpha_\lambda + \alpha_\lambda \wedge \alpha_\lambda = 0
\textup{ for all $\lambda$}   \hspace{.3cm} \Leftrightarrow   \hspace{.3cm}
\textup{(\ref{eq1}), (\ref{eq2}) and (\ref{eq3}) plus flatness.}
\edm
For each real value of $\lambda$ we can integrate $\alpha_\lambda$ to obtain a frame
for a flat immersion.
Thus the flatness condition can be encoded by assuming that we have such a 1-parameter
\emph{family} of frames.\\

In general, let $G$ be a complex semisimple Lie group, and suppose we have the following ingredients:
\begin{enumerate}
  \item for $\lambda \in \C^*$, a 1-parameter {family} of 1-forms, 
    $\alpha_\lambda \, \in \mathfrak{g} \otimes \Omega (M)$. 
  \item  $\alpha_\lambda$ is a Laurent polynomial in $\lambda$,
  \bdm 
    \alpha_\lambda = \sum_{i=a}^b a_i \lambda^i,  
          \hspace{1cm} a_i \in \mathfrak{g} \otimes \Omega (M).
  \edm
  \item 
  $\alpha_\lambda$ satisfies the Maurer-Cartan equation for all $\lambda \in \C^*$.
\end{enumerate} 
Then we can integrate to obtain family $F_\lambda: M \to G$, and project
 to obtain a family of special submanifolds $f_\lambda: M \to G/H$, where $H$ is
 some subgroup of $G$.\\
   \textbf{Interesting question:}  what are the special submanifolds corresponding to
   the projections $f_\lambda$?

%!!!!!!!!!!!!!!!!!!!!!!!!!!!!!!!!!!!!!!!!!!!!!!!!!!!!!!!!!!!!!!!!!!!!!!!
%!!!!!!!!!!!!!!!!!!!!!!!!!!!!!!!!!!!!!!!!!!!!!!!!!!!!!!!!!!!!!!!!!!!!!!!!
\begin{figure}[here]
\includegraphics[width=12cm]{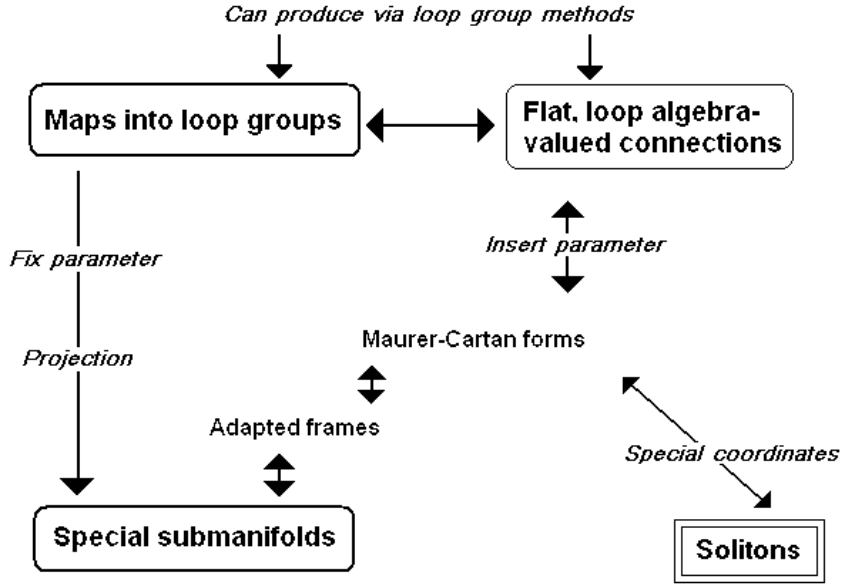}\caption{The relations between maps into
loop groups, flat connections, special submanifolds and special PDE.}
\end{figure}

\subsection{The Connection to Special PDE}
The existence of a 1-parameter family of integrable Maurer-Cartan forms
(corresponding to flat connections with values in a loop algebra)
is well known to be an essential characteristic of soliton equations
and other so-called integrable systems. This aspect manifests itself in the
following way: given a family of 1-forms
$\alpha_\lambda = \sum_{i=a}^b a_i \lambda^i, ~ 
           a_i \in \mathfrak{g} \otimes \Omega (M)$, 
as above, it is easy to see that
\bdm 
\dd \alpha_\lambda + \alpha_\lambda \wedge \alpha_\lambda = 0, \hspace{.5cm} \textup{for all } \lambda
\edm
if and only if
\bdm 
\dd a_k + \sum_{i+j = k} a_i \wedge a_j = 0.
\edm
This is a system of PDE (after choosing some coordinates).

\textbf{Example:} We return once more to our example of flat immersions into $S^3$.
The Gauss equation:
$\dd \omega + \omega \wedge \omega -\beta \wedge \beta^t -\theta \wedge \theta ^t = 0$,
together with the flatness condition $\dd \omega + \omega \wedge \omega = 0$,
turn out to reduce to one equation, in special coordinates:
\bdm
\phi_{xy} = 0, 
\edm
namely, the wave equation.

%*********************************************
\section{Grassmann Geometries}
The methods from loop groups used here produce submanifolds which are, or are
related to, Grassmann submanifolds in homogeneous spaces.  This point has perhaps
not been emphasized in the past, because the majority of applications studied 
were in space forms, where the Grassmann 
submanifold condition (arising from orbits of the action of the isometry group 
in the symmetric space representation) is satisfied by any submanifold.

The concept of a Grassmann submanifold was introduced by Harvey and Lawson in \cite{harveylawson}, as follows:
let $\NN$ be a manifold and take 
any subset, $\VV$,  of the
Grassmann bundle over $\NN$ consisting of tangential $s$-planes, 
$Gr_s(T\NN) = \cup_{x \in \NN} Gr_s(T_x\NN)$.  
A \emph{$\VV$-submanifold}, $N$, of $\NN$, is an $s$-dimensional connected 
submanifold such that $T_xN \in \VV$ for each $x \in N$.
The set of such submanifolds, $N$, is called the  \emph{$\VV$-geometry}.  

In this article,  $\NN$ will always be a homogeneous space, $G/H$, with 
$G$ a connected Lie group, and
 $\VV$ an orbit of the action of $G$ on  $Gr_s(T\NN)$. In such a case, the
geometry $\VV$ is determined by an $s$-dimensional vector subspace 
of the tangent space at the origin, $H$, of of $G/H$. A special case is
when
 $\NN = \Ubar/\Kbar$ is a \emph{symmetric space}. Then we have
 the  canonical decomposition of the Lie algebra 
$\ubar = \kbar \oplus \pbar$, and the tangent space at the origin is
  $T_0 \NN =  \pbar$. So for symmetric spaces we have the correspondence:
\bdm
\{\textup{$s$-Dim $\VV$-geometries} \} \leftrightarrow \{\textup{$s$-Dim subspaces } \pp \subset \pbar \}.
\edm
Given $\pp \subset \pbar$, we will call  the associated geometry the $\VV_\pp$-geometry.

If $Ad_{\Kbar} \pp \subset \pp$ then the $\VV_\pp$-geometry consists of integral
submanifolds of a \emph{distribution} determined by $\pp$, but otherwise it is
a more general concept.

\subsection{Examples}
For space forms, \emph{any} $s$-dimensional submanifold is a $\VV_\pp$-submanifold
for any $s$-dim subspace $\pp \subset \pbar$. We demonstrate this for curves in
 $\NN = SO(3)/SO(2) = S^2$. We have the canonical decomposition:
\beqas
\mathfrak{so}(3) = \kbar \oplus \pbar = \left\{ \bbar {ccc} * & * & 0 \\ * & * & 0\\0 & 0 & 0 \ebar \right\} \oplus
  \left\{ \bbar {ccc} 0 & 0 & * \\ 0 & 0 & *\\ {*} & * & 0\\ \ebar \right\}.
\eeqas
For a 1-dimensional subspace $\pp \subset \pbar$, we can, with no loss of generality,
take $\pp = \left\{ \bbar {ccc} 0 & 0 & * \\ 0 & 0 & 0\\ {*} & 0 & 0\\ \ebar \right\}$.  

Let $f: \real \to S^2$ be any curve.
The $\VV_\pp$-geometry is determined by the left action of $SO(3)$ on
$Gr_1(T S^2)$, and to show that a curve in $S^3$ is a $\VV_\pp$-submanifold,
we need to show there exists 
 frame $F \in SO(3)$ for $f$, such that the projection onto $\pbar$ of $F^{-1} \dd F$
lies in  $\pp$.
This is achieved by choosing an \emph{adapted frame} $F: \real \to SO(3)$,
\bdm
F = [e,n,f], \hspace{1cm} \textup{$e$ tangent, $n$ normal},
\edm
\bdm
F^{-1} \dd F = \bbar{c} e^t\\n^t\\f^t \ebar \bbar {ccc}\dd e & \dd n & \dd f \ebar =
\bbar {ccc} 0 & e^t\dd n & e^t \dd f \\ n^t \dd e & 0 & n^t \dd f \\ f^t \dd e & f^t \dd n & 0 \ebar.\\
\edm
The $\pbar$ part is $\bbar {ccc} 0 & 0 & e^t \dd f \\ 0 & 0 & n^t \dd f \\ f^t \dd e & f^t \dd n & 0 \ebar = \bbar{ccc} 0 & 0 & e^t \dd f \\ 0 & 0 & 0 \\ f^t \dd e & 0 & 0 \ebar \in \pp$.\\

More meaningful examples of Grassmann submanifolds are
Lagrangian submanifolds of $\cc P^n$ and
almost complex and totally real submanifolds of $S^6$. The latter arise with
respect to the action of $G_2$ on the homogeneous space $S^6 = G_2/SU(3)$,
which is not a symmetric space representation of $S^6$; hence there is no
conflict with the above comment concerning space forms.

\section{Grassmann geometries associated to loop groups}
Loop group techniques (AKS-theory, DPW, etc) produce maps into a subgroup of a
loop group which are characterized by the fact that the Maurer-Cartan form is
a Laurent polynomial of fixed degree in the loop parameter, $\lambda$.
Solutions are determined modulo the action of the constant
subgroup - hence actually frames for maps into a homogeneous space.

We formulate this in the language of Grassmann geometries:
Let $G$ be a complex semisimple Lie group, and define the
loop group 
\bdm
\Lambda G := \{ \gamma: S^1 \to G \},
\edm 
where the maps have some convergence condition, such as the Wiener topology,
which makes $\Lambda G$ a Banach Lie group.
Let $\hh$ be a Banach subgroup of $\Lambda G$, and denote by 
 $\hh^0 := \hh \cap G$, the subgroup of constant loops.
Then the left coset space $\hh/\hh^0$ is a  
homogeneous space on which $\hh$ acts on the left.

To define Grassmann geometries on $\hh/\hh^0$, we need to describe its tangent space at the
origin. The Lie algebra of $\Lambda G$ is 
$\Lambda \mathfrak{g} = \{\sum_{i-\infty}^\infty a_i \lambda^i ~|~ a_i \in \mathfrak{g}\}$,
 and $\textup{Lie}(\hh)$ is a vector subspace of  $\Lambda \mathfrak{g}$.
Clearly
 $\textup{Lie}(\hh^0) = \{$ constant polynomials in $\textup{Lie}(\hh) \}$, from
 which it follows that
\bdm
T_0 \frac{\hh}{\hh^0} = \{ \sum_{i \neq 0} a_i \lambda^i \} \subset \textup{Lie}(\hh).
\edm

For integers $a<b$, define
 $W_a^b \subset T_0 \frac{\hh}{\hh^0}$  by
\bdm
W_a^b = \{ x \in T_0 \frac{\hh}{\hh^0} ~|~ \sum_{i=a}^b a_i \lambda^i \}.
\edm
Now set $\VV_a^b$ to be the distribution given by the 
orbit of $W_a^b$ under the action of $\hh$
on $Gr_{b-a}\big (T \frac{\hh}{\hh^0} \big )$.

The basic object we can construct, using  the techniques described here,
are $\VV_a^b$-compatible
(immersed) submanifolds of $\hh/\hh^0$,
i.e. maps $f: M \to \hh/\hh^0$ for which there exists  frames 
$F: M \to \hh$ with $F^{-1} \dd F = \sum_{i=a}^b \alpha_i \lambda^i$.

\section{Special submanifolds from loop group maps}

A $\VV_a^b$-immersion $f: M \to \hh/\hh^0$, leads naturally to families of special
submanifolds as follows:
Evaluate $f$ at some $\lambda_0$, to get a map $f_{\lambda_0}: M \to G/\hh^0$.
The subgroup $\hh$ together with the $\VV_a^b$ condition make $f_{\lambda_0}$ a
certain Grassman submanifold.

Since $f$ is a $\VV_a^b$-immersion, by definition, there exists a lift
 $F: M \to \hh$, such that
 $\alpha := F^{-1} \dd F = \sum_{i=a}^b \alpha_i \lambda^i$.
An essential point is: $\alpha$ must satisfy the Maurer-Cartan equation,
\bdm
\dd \alpha + \alpha \wedge \alpha = 0,
\edm
 for \emph{all}
values of $\lambda$.
This is equivalent to some conditions on $\alpha_i$, 
\beq \label{mccondition}
\dd \alpha_k + \sum_{i+j = k} \alpha_i \wedge \alpha_j = 0,
\eeq
independent of $\lambda$.

The equations (\ref{mccondition}) give some \emph{extra conditions}, usually on the
(tangent and normal) curvature of the submanifold. This will be illustrated by our
example below.

%**************************************************************
\section{The three involution loop group} \label{construction}
Now we define the generalization of the loop group construction
of \cite{feruspedit1996}.
Let $G$  be a complex semisimple Lie group and $\tbar$, $\shat$, $\rho$ 
 commuting involutions of $G$, where $\rho$ is $\cc$-antilinear.
The fixed point subgroup with respect to $\rho$, $\Ubar := G_\rho$, is a real form of
the group.

We extend the involutions to $\Lambda G$ by the rules:
\beqas \label{symrho}
(\rho X)(\lambda) = \rho(X(\bar{\lambda})),\\
(\shat X)(\lambda) = \shat(X(-\lambda)), \label{symshat}\\
(\tbar X)(\lambda) = \tbar(X(-1/\lambda)), \label{symtau}
\eeqas
and consider the subgroup fixed by all three involutions:
\bdm
\hh = \Lambda G_{\rho \tbar \shat}.
\edm

Consider a $\VV_{-1}^1$-immersion $f:M \to \hh/\hh^0$.
For $\lambda \in \real^*$, 
\bdm
f_\lambda :M \to \Ubar/\Ubar_{\tbar \shat},
\edm  since
$\hh^0 = \Ubar_{\tbar \shat} =\Ubar_\tbar \cap \Ubar_{\shat}$. \\

We can also project to obtain 
maps into the symmetric spaces $\Ubar/\Ubar_{\tbar}$ and $\Ubar/\Ubar_{\shat}$, or more
generally, into any homogeneous space $\Ubar/H$, where $\Ubar_\tbar \cap \Ubar_{\shat} \subset H$. \\

 What are the special submanifolds so obtained?

\section{Reflective submanifolds}
We are primarily interested in the projection to $\Ubar/\Ubar_{\tbar}$, as this
generalizes the isometric immersions of space forms studied in \cite{feruspedit1996}.
To describe the projections, we first need to define reflective submanifolds.\\

\textbf{Examples} In space forms, these are just the complete totally geodesic submanifolds.
 Other examples are Lagrangian embeddings of 
$\real P^n \subset \cc P^n$ and $\real H^n \subset \cc H^n$. \\

 \textbf{Definition:}
A \emph{reflective submanifold}, $N$, of a Riemannian manifold, $\bar{N}$,
is a totally geodesic symmetric submanifold.\\

For a connected symmetric space
$\bar{N} = \bar{U}/\bar{K}$, we can characterize a reflective submanifold $N$ of $\bar{N}$,
by the existence of a second involution on the Lie algebra of $\bar{U}$.  
Specifically, $N \subset \bar{N}$ is  
characterized by a pair of commuting involutions, $\bar{\mathfrak{\tau}}$ and $\hat{\sigma}$,
of the Lie algebra $\bar{\mathfrak{u}}$ of $\bar{U}$, and $\bar{K} = \bar{U}_{\bar{\mathfrak{\tau}}}$.
 That is:
\bdm
N \subset  \bar{U}/\bar{K}\hspace{.5cm} \leftrightarrow \hspace{.5cm} 
(\bar{\mathfrak{u}}, \bar{\mathfrak{\tau}}, \hat{\sigma}).
\edm
We have two  canonical decompositions of the Lie algebra
 $\bar{\mathfrak{u}} = \bar{\mathfrak{k}} \oplus \bar{\mathfrak{p}} = \hat{\mathfrak{k}} \oplus \hat{\mathfrak{p}}$, into the $+1$ and $-1$ eigenspaces of the two involutions.
  Setting 
 \bdm
 \mathfrak{p} := \bar{\mathfrak{p}} \cap \hat{\mathfrak{p}},
\edm
  the reflective submanifold is given by:  
$N = \pi_{\bar{N}} \exp(\mathfrak{p})$. \\

Reflective submanifolds of symmetric spaces were classified by DSP Leung
(1974-1979), and there are clearly many cases.

\section{{Isometric immersions of space forms}}
The  three involution loop group leads naturally to  a 
 generalization of the following results/conjectures:\\

An isometric immersion $f: M^k(c) \to M^n(\tilde{c})$, of space forms with 
constant sectional curvature $c$ and $\tilde{c}$ respectively,
 has \emph{negative extrinsic curvature} if $c< \tilde{c}$. 
  There are two basic questions: existence of a local solution,
and existence of a complete solution. For these it is known:

\begin{enumerate}
\item Local solutions exist iff $n \geq 2k-1$ (Cartan).\\
\item 
Theorem (JD Moore): If $0<c<1$, there is no complete isometric immersion with flat normal
bundle of $S^k(c)$ into $S^n$ for any $k>1$ and any $n$. \\
\item 
Plausible conjecture: If $c<-1$ there is no complete isometric immersion with flat 
normal bundle of $H^k(c)$ into $H^n(-1)$ for any $k>1$ and any $n$.

For the case $n = 2k-1$, this is equivalent to the conjectured generalization
of  Hilberts's non-immersibility of $H^2$ into $E^3$.
\end{enumerate}

\section{The generalization to other reflective submanifolds}
$M$ a Riemannian manifold, let $M_R$ denote the same manifold
with the metric scaled by a factor $R>0$. \\

\noindent \textbf{Problem A:} 
Suppose given a reflective submanifold
\bdm
N \subset \bar{N}
\edm
 of a symmetric space.
  Thus, $N_R \subset \bar{N}_R$ is also a reflective submanifold.
   Does there exist a (local or global) isometric immersion 
\bdm
N_R \to \bar{N},
\edm
 satisfying condition X? That is, can we shrink/stretch $N$ within $\bar{N}$?
More specifically, we ask this for: 
\begin{enumerate}
\item
$R>1$, if $\bar{N}$ is of compact type,
\item
$R<1$, if $\bar{N}$ is of non-compact type. \\
\end{enumerate}

Reflective submanifolds in other symmetric spaces do not
generally have flat normal bundle. Thus, 
we need to replace the flat normal bundle condition with an appropriate
one, which we call here \emph{condition X}. \\

Condition X just says: 
\begin{enumerate}
\item $N_R \to \bar{N}$ is a $\VV_\pp$-submanifold, where $N = \exp(\pp)$.\\
\item
The normal bundle of $N_R \to \bar{N}$ is 
isomorphic (as an affine vector bundle/connection pair) with the normal
bundle of $N_R \subset \bar{N}_R$. 
\end{enumerate}

\section{Projections to $\Ubar/\Ubar_{\tbar}$ and $\Ubar/\Ubar_{\shat}$}
%********************************
Here we summarize results from \cite{reflective}. In fact Proposition \ref{prop2}
is stated incorrectly in \cite{reflective} - the limit as $\lambda \to \infty$ or
$\lambda \to 0$ must be taken before a curved flat is obtained.

Set $\Kbar := \Ubar_{\tbar}$, and $\Khat := \Ubar_{\shat}$.
Take $f: M \to \mathcal{H}/\mathcal{H}^0$ a $\mathcal{V}_{-1}^1$-immersion. 
Recall $f_\lambda: M \to \Ubar/(\Kbar \cap \Khat)$, for $\lambda \in \real$.

\begin{proposition} 
Let $\bar{f}_\lambda: M \to \bar{U}/\Kbar$
be the projection of $f_\lambda$.   Suppose that
 $\bar{f}_\lambda$ is regular. Then $\bar{f}_\lambda$ is a solution of Problem A 
 (for $R>1$). Conversely, any solution of Problem A, corresponds to
such a $\mathcal{V}_{-1}^1$-immersion.
\end{proposition}
%*****************************
\begin{proposition}  \label{prop2} 
Let $\hat{f}_\lambda: M \to \bar{U}/\Khat$
be the projection of $f_\lambda$. Then:
\begin{itemize} 
\item $\hat{f}_\lambda$ is asymptotic to a curved flat
in $\bar{U}/\Khat$, as $\lambda \to \infty$, or $\lambda \to 0$.
\item
If $\bar{f}_\lambda$ is regular then so is $\hat{f}_\lambda$ (but not conversely).\\

Hence, if $\Ubar/\Kbar$  \textbf{compact} then:
\begin{enumerate}
\item Local regular solutions to Problem A exist\\$\Rightarrow$ 
$\textup{Dim}(\mathfrak{p}) \leq \textup{Rank}(\bar{U}/\Khat)$. 
\item Global regular solutions to Problem A do not exist 
for Dim($N) > 1$.
\end{enumerate}
\end{itemize}
\end{proposition}

\section{Consequences}
\begin{theorem} \label{compactthmbrander}
\textbf{(Compact Case)} The following list contains the geometric interpretations of all possible solutions
 to Problem A
for the case $R>1$ and  $\bar{N}$ is a simply connected, compact, irreducible, 
Riemannian symmetric space. 

 In all cases, local solutions exist and can be constructed by
loop group methods. In all cases where $\textup{Dim}(N_R) >1$,
 there is no solution which is geodesically
complete.
\begin{enumerate}
\item
$N_R=S_R^k$ is an isometric immersion with flat normal bundle of a $k$-sphere of radius
$\sqrt{R}$ into the unit sphere $S^n$, with $0< k \leq (n+1)/2$, and $n \geq 2$.
\item
$N_R=S_R^n$ is an  isometric totally real immersion of an $n$-sphere of radius $\sqrt{R}$
into complex projective space ${\bf C} P^n$, with $n \geq 2$.
\end{enumerate}
\end{theorem}
Note: Lagrangian immersions of a sphere
into ${\bf C} P^n$ is a new example of a submanifold as 
an integrable system.\\

\begin{theorem}
\textbf{(Non-Compact Case)}
The analogue - \textbf{except} we do not obtain the global non-existence result,
which remains an open problem.
\end{theorem}

\nopagebreak
%*********************************

\end{document}